\documentclass[11pt]{article}

\usepackage{amssymb}
\usepackage{amssymb}
\usepackage{epsfig}
\usepackage{amsmath}
\usepackage{color}
\usepackage[english]{babel}
\usepackage[latin1]{inputenc}

\newtheorem{theorem}{Theorem}[section]
\newtheorem{lemma}[theorem]{Lemma}
\newtheorem{e-proposition}[theorem]{Proposition}
\newtheorem{corollary}[theorem]{Corollary}
\newtheorem{e-definition}[theorem]{Definition\rm}
\newtheorem{remark}{\it Remark\/}
\newtheorem{example}{\it Example\/}

\begin{document}
\centerline{\bf Microlocal Versal Deformations of Plane Curves}
\bigskip
\centerline{Jo\~ao Cabral and Orlando Neto\footnote{This research was partially suported by \em FEDER \em and \em FCT-Plurianual 2009 \em.}}
\bigskip

\begin{abstract}
We introduce the notion of microlocal versal deformation of a plane curve. 
We construct equisingular versal deformations of Legendrian curves that are the conormal of a semi-quasi-homogeneous branch.
\end{abstract}

\section{Contact Geometry}
\noindent
Schlessinger, Tyurina and Grauert initiated the study of versal deformations of analytic spaces 
(cf. \cite{GP5}).
We show in remark \ref{final} that the obvious definition of deformation of a 
Legendrian curve is not very interesting.
There would be too many rigid Legendrian curves. 
We consider an alternative approach, recovering Lie's  original point of view: 
to look at contact transformations as maps that take plane curves into plane curves.

\noindent
We follow the terminology of \cite{GP5} regarding deformations of analytic spaces. 
When we refer to a curve or a deformation, we are identifying it with its germ at a convenient point. 
If $Y\hookrightarrow \mathcal Y \to T$ is a deformation $\mathcal Y$ of a plane curve $Y$ over a complex manifold $T$,
we assume that $Y$ is reduced. 
We say that a deformation of a plane curve is \em equisingular \em if all of its fibers have the same topological type.
Given a germ of an analytic set $Y$ we will denote by mult $Y$ the multiplicity of $Y$.

\noindent
We set $\mathbb C^A=\{(t_a)_{a\in A}: t_a\in \mathbb C\}$, $t_A=(t_a)_{a\in A}$, for each finite set $A$.
In general, we denote $\partial f/\partial t$ by $\partial_t f$.


\noindent
Let $(X,\mathcal O_X)$ be a complex manifold. 
A differential form of degree $1$ is called a \em contact form
\em if $\omega\wedge d\omega$ never vanishes.
If $\omega$ is a contact form there is a  system of local coordinates $(x,y,p)$ such that
$\omega = dy-pdx$. 
A locally free $\mathcal O_X$-module  $\Omega$ of rank $1$ is called a \em contact structure \em 
if $\Omega$ is locally generated by a contact form. The pair $(X,\Omega)$ is called a \em contact threefold. \em
\noindent
Let $\xi dx+\eta dy$ denote the canonical $1$-form of $T^*\mathbb C^2=\mathbb C^2\times \mathbb C^2$. 
We will identify the open set $\{\eta\not=0\}$
of $\mathbb P^*\mathbb C^2=\mathbb C^2\times \mathbb P^1$ with the vector space $\mathbb C^3$, 
endowed with the coordinates $(x,y,p)$, where $p=-\xi/\eta$. 
The contact structure $\Omega$ of $\mathbb C^3$ generated by the differential 
form $\omega=dy-pdx$ is the restriction to $\mathbb C^3$ of the canonical contact structure of $\mathbb P^*\mathbb C^2$. 
Let $\pi:\mathbb P^*\mathbb C^2\to \mathbb C^2$ be the canonical projection.
Let $C(Z)$ be the tangent cone of a germ of an analytic set $Z$.
A Legendrian curve $L\subset (\mathbb P^*\mathbb C^2,\sigma)$ is in \em strong generic position \em if $C(L)\cap (D\pi(\sigma))^{-1}(0)=\{ 0\}$.

\noindent
\begin{remark}\label{generalidades}\em
\begin{enumerate}
\item[$(i)$]
Let $Y_1$, $Y_2$ be germs of irreducible plane curves at $o$. Let $(L_i,\sigma_i)$ be the conormal of $Y_i$, $i=1,2$. 
Then $\sigma_1=\sigma_2$ if and only if $Y_1$ and $Y_2$ have the same tangent cone. 
Hence the plane curve $\pi (L)$ has irreducible tangent cone for each Legendrian curve $L\subset \mathbb P^*\mathbb C^2$.
We will assume that all plane curves have irreducible tangent cone.
\item[$(ii)$]
If $Y\hookrightarrow \mathcal Y \to T$ is an equisingular deformation of a plane curve $Y$, 
$Y\hookrightarrow \mathcal Y \to T$ is isomorphic to a deformation $Y\hookrightarrow \mathcal Z \to T$  
such that the tangent cone of $\mathcal Z_{t}$ does not depend on $t$. 
We will assume that all the fibers of an equisingular deformation of a plane curve have the same tangent cone.
\item[$(iii)$]
Let $Y$ be a germ of an irreducible plane curve with tangent cone $\{y=0\}$. 
Let $y=\sum_{\alpha\in\mathbb Q} a_{\alpha} x^{\alpha}$ be the Puiseux expansion of $Y$. 
Set $\delta=\inf \{\alpha:a_{\alpha}\neq 0\}$. The tangent cone of the conormal of $Y$ equals 
$\{x=y=0\}$ if $\delta<2$, $\{y=p-2a_2 x=0\}$ if $\delta=2$ and $\{y=p=0\}$ if $\delta>2$.
\end{enumerate}
\end{remark}

\noindent
Let $(X,\Omega_X)$ be a contact threefold. Let $T$ be a complex manifold. 
Let $q_X:X\times T\rightarrow X$, $r_X:X\times T\rightarrow T$ be the canonical projections. 
The pair $(X\times T,q_X^{\ast}\Omega_X)$ is called a \em relative contact threefold. \em 
Let $\mathcal{L}$ be an analytic set of dimension  $dim\: T+1$ of $X\times T$. 
We say that $\mathcal{L}$ is a  \em  relative Legendrian curve \em if 
$q_X^{\ast}\omega$ vanishes on the regular part of  $\mathcal{L}$, for each section $\omega$ of $\Omega_X$. 
The intersection of  $\mathcal{L}$ with each fiber of $r_X$ is a Legendrian curve.
 Given two relative contact threefolds $X\times T$ and $Y\times T$, a biholomorphic map 
$\Phi:X\times T\rightarrow Y\times T$   is called a 
\em relative contact transformation \em if $r_Y \circ \Phi=r_X$ and $\Phi^*q_Y^*\Omega_Y=q^*_X\Omega_X$. 
Given $t\in T$, let $\Phi_t$ be the induced map from $X\times \{t\}$ into $Y\times \{t\}$.
Two relative Legendrian curves $\mathcal L_1\subset X\times T,\mathcal L_2\subset Y\times T$ are \em isomorphic \em
if there is a relative contact transformation $\Phi: X\times T\to Y\times T$ such that $\Phi(\mathcal L_1)=\mathcal L_2$.
Let $Y\hookrightarrow \mathcal Y \to T$ be a deformation of a plane curve. 
Let $\pi_T:\mathbb P^*\mathbb C^2\times T\to \mathbb C^2\times T$ be the canonical projection.
The \em conormal \em of $ \mathcal{Y}$ is the smallest relative Legendrian curve $\mathcal L\subset \mathbb P^*\mathbb C^2\times T$ such that
$\pi_T(\mathcal L)=\mathcal Y$.

\noindent
Let $\mathcal{H}$ be the group of contact transformations of the type $(x,y,p)\mapsto (\lambda x, \mu y, \mu \lambda^{-1} p)$, $\lambda,\mu\in \mathbb C^{\ast}$. Set $\ell=\{y=p=0\}$. Let $\mathcal{G}$ be the group of germs of contact transformations  $\varphi :(\mathbb C^3,0)\to (\mathbb C^3,0)$ such that $D\varphi(0)(\ell)=\ell$. Let $\mathcal{J}$ be the group of germs of contact transformations 
\begin{equation}\label{grupojota}
(x,y,p)\mapsto(x+\alpha,y+\beta,p+\gamma), \ \hbox{\rm where} \  
\alpha,\beta,\gamma,\partial_x \alpha,
\partial_y \beta,\partial_p \gamma\in(x,y,p).
\end{equation}

\begin{theorem}\label{inv}
\em (cf. \cite{ARNE}) \em The group $\mathcal{J}$ is an invariant subgroup of $\mathcal{G}$. 
Moreover, $\mathcal{G}/\mathcal{J}$ is isomorphic to $\mathcal{H}$.
\end{theorem}

\begin{theorem}\label{ABC}
\em (cf. \cite{ARNE}) \em Set $t=(t_1,\ldots,t_m)$. Let $\alpha\in\mathbb{C}\{x,y,p,t\}$ and $\beta_0\in\mathbb{C}\{x,y,t\}$ 
be power series such that $\alpha,\partial_x \alpha,\beta_0,\partial_y \beta_0\in (x,y,p)$. 
There are $\beta,\gamma\in \mathbb{C}\{x,y,p,t\}$ such that $\beta-\beta_0\in (p)$, $\gamma\in (x,y,p)$ and (\ref{grupojota}) is a relative contact transformation. Moreover, $\beta$ is the solution of the Cauchy problem
\begin{equation}\label{cauchy}
\left(1+\frac{\partial\alpha}{\partial x}+p\frac{\partial\alpha}{\partial y}\right) \frac{\partial\beta}{\partial p}-p\frac{\partial\alpha}{\partial p}\frac{\partial\beta}{\partial y}-\frac{\partial\alpha}{\partial p}\frac{\partial\beta}{\partial x}=p\frac{\partial\alpha}{\partial p}, ~~~~ \beta-\beta_0\in (p). 
\end{equation}
\end{theorem} 

\noindent
We set $\alpha=\sum_{l\geq 0} \alpha_l p^l$ and $\beta=\sum_{l\geq 0} \beta_l p^l$, where $\alpha_l,\beta_l\in\mathbb C\{x,y,t\}$ for all $l$.

\noindent
\begin{theorem}\label{legequi}
Let $L_1,L_2$ be two germs of Legendrian curves of $\mathbb P^{\ast}\mathbb{C}^{2}$ 
in strong generic position. Set $Y_i=\pi(L_i)$, $i=1,2$. 
If there is a contact transformation $\varphi$ such that $\varphi(L_1)=L_2$, 
there is a local homeomorphism $\psi:(\mathbb C^2,0)\rightarrow (\mathbb C^2,0)$ such that $\psi(Y_1)=Y_2$.
\end{theorem}

\begin{proof}
Assume that $L_1$ is irreducible. Assume that $(n_i,k_i)=1$, $1\le i\le l$, and $n_1/k_1<\cdots<n_l/k_l$. 
Set $k_0=1,k=k_l$.
The curve $Y_1$ has Puiseux pairs $n_i/k_i$, $1\le i\le l$, if and only if the following conditions hold for $1\le i\le l$: 
$($i$)$ $a_{n_ik/k_i}\not=0$; 
$($ii$)$ if $a_j\not=0$ and $k/k_{i-1}$ does not divide $j$, $j\ge n_ik/k_i$.

\noindent 
 Let $i\in\{1,2\}$. 
Following the notations of remark 1(iii), we can assume that $C(Y_i)=\{y=0\}$, 
$\delta> 2$ and the tangent cone of $L_i$ equals $\{y=p=0\}$. 
By theorem \ref{inv},  we can assume
that $\varphi$ is of the type (\ref{grupojota}). Composing $\varphi$ with the contact transformation induced by $(x,y)\mapsto (x,y-\beta_0)$, 
we can assume that $\beta_0=0$. Hence $\partial_x \beta,\partial_y \beta \in (p)$. 
Setting $p=0$ in (\ref{cauchy}), we conclude that $\beta_1=0$. Let
\begin{equation}\label{expfamilia2}
x=s^{k}, \qquad y=s^{n}+{ \sum_{r\geq n+1}a_r s^{r}},\qquad p=\frac{n}{k}s^{n-k}+{ \sum_{r\geq n+1}}\frac{r}{k}a_r s^{r-k}.
\end{equation}
be a parametrization of $L_1$. 
The curve $Y_2$ admits a parametrization 
 of the type $x=s^k+\alpha(s),y=s^n+\sum_{r\geq n+1} a_r s^r+\sum_{l\geq 2}\beta_l(s) p(s)^l=s^k+\sum_{r\geq n+1} b_r s^r$, where $\alpha\in (s^{k+1})$. Since $2(n_i-k_i)>n_i$, the coefficients of $\beta_l(s)p(s)^l$, $l\ge 2$, verify condition $($ii$)$. Hence the coefficients $b_r$ verify conditions $($i$)$,$($ii$)$.
Replacing the parameter $s$ by a parameter $t$ such that $t^k=s^k+\alpha(s)$ 
one obtains a parametrization of $Y_2$ that still verifies conditions $($i$)$ and $($ii$)$ (cf. Lemma 3.5.4 of \cite{CTCW}).

\noindent
If $L'_1$ and $L''_1$ are irreducible Legendrian curves, a similar argument shows that the contact order (cf. Section 4.1 of \cite{CTCW}) of $\pi(\varphi(L'_1))$ 
and $\pi(\varphi(L''_1))$ equals the contact order of $\pi(L'_1)$ and $\pi(L''_1)$ .$\square$
\end{proof}

\begin{e-definition}\label{planeproj}\em
Let $\mathcal L$ be a relative Legendrian curve of a relative contact threefold $X\times (T,o)$. 
We call  \em plane projection \em of  $\mathcal L$ 
to an analytic set $\pi_T(\Phi(\mathcal L))$, 
where $\Phi: X\times T\to \mathbb P^{\ast}\mathbb{C}^2\times T$ is a relative contact transformation. 
The plane projection $\pi_T (\Phi (\mathcal L))$ is called \em generic \em if 
$\Phi_o(\mathcal L_o)$ is in strong generic position.
Two Legendrian curves are  \em equisingular \em if their generic plane projections have the same topological type.
\end{e-definition}

\begin{lemma}\label{minimal}
If $Y$ is a plane projection of a Legendrian curve $L$, mult $Y\ge$ mult $L$.
Moreover,   mult $Y=$ mult $L$ if and only if $Y$ is a generic plane projection of $L$.
Hence the multiplicity is an equisingularity invariant of a Legendrian curve.
\end{lemma}

\begin{proof}
Its enough to prove the result when $L$ is irreducible.
Let $\varphi :X\to \mathbb P^*\mathbb C^2$ be a contact transformation.
Let (\ref{expfamilia2}) be a parametrization of $\varphi (L)$.
Then mult $Y=k$. Moreover, mult $L=k$ if $n\ge 2k$ and mult $L=n-k$ if $n< 2k$.
The result follows from remark \ref{generalidades}(iii).
$\square$
\end{proof}

\section{Microlocal Versal Deformations}

\noindent
Let $L$ be a Legendrian curve. A relative Legendrian curve $\mathcal L\subset X\times (T,o)$ is
called a \em  deformation \em of $L$ if the sets $\mathcal L_o$ and $L$ are equal. 
The deformation $\mathcal L$ is called \em equisingular \em if its fibers are equisingular.
We do not demand the flatness of the morphism $\mathcal L\hookrightarrow X\times T\to T $.
Let $Y$ be a plane curve with conormal $L$. A deformation $\mathcal Y$ of $Y$ is called a \em microlocal deformation \em of
$Y$ if the conormal of $\mathcal Y$ is a deformation of $L$.

\noindent
Let $\mathcal J'_T$ be the group of relative contact transformations
 $\Phi: (\mathbb C^3,0)\times (T,o)\to (\mathbb C^3,0)\times (T,o)$ such that $\Phi_o=$id$_{\mathbb C^3}$.
 Let $m$ be the maximal ideal of $\mathcal O_{T,o}$. Given $\Phi\in \mathcal J'_{T}$, 
there are $\alpha,\beta,\gamma\in m$ such that $\Phi$ equals (\ref{grupojota}).
 Two deformations $\mathcal L_1,\mathcal L_2$ of $L$ over $(T,o)$ are \em isomorphic \em
if there is  $\Phi \in \mathcal J'_{T}$ such that $\Phi (\mathcal L_1)=\mathcal L_2$.
Two microlocal deformations of a plane curve are \em microlocally equivalent \em if 
their conormals  are isomorphic.

\begin{remark}\label{final}\em
Let $\mathcal L$ be a flat equisingular deformation of the conormal of  $\{y^k=x^n\}$ along $(T,o)$, where $n>k>1$ and $(k,n)=1$. 
Since $\mathcal L$ is flat and $ny-kxp$ vanishes on $\mathcal L_o$, 
there is an holomorphic function $f$ on $X\times T$ such that $f$ vanishes on $\mathcal L$ and $f(x,y,p,o)=ny-kxp$.
Hence, if $t\in T$ is close enough to $o$, $\mathcal L_t$ is contained in a smooth hypersurface. By theorem 8.3 of \cite{SKKO}, 
there are integers $l,m$ and a  system of local coordinates $(x',y',p')$ on $X$ such that 
$m>l>1$, $(l,m)=1$ and $\mathcal L_t$ is the conormal of $\{y'^l=x'^m\}$.
Since the deformation is equisingular, $l=k$ and $m=n$. Hence $\mathcal L_t$ is isomorphic to $\mathcal L_o$ when $t$ is close enough to $o$.
Modifying the proof of the theorem refered above one can show that de deformation $\mathcal L$ is trivial even without the assumption that $\mathcal L$ is equisingular.
\end{remark}


\begin{theorem}\label{lematwo}
Let $\mathcal L$ be an equisingular deformation of a Legendrian curve $L$ over $(T,o)$. 
A generic plane projection of $\mathcal L$ is an equisingular deformation of $\pi(L)$.
\end{theorem}
\begin{proof}
Let $\Phi:X\times T\to \mathbb P^*\mathbb C^2\times T$ be relative contact transformation 
such that $\Phi_o(\mathcal L_o)$ is in strong generic position.
It follows from lemma \ref{minimal} and the upper-semicontinuity of the map 
$t\mapsto$ mult $\pi(\Phi_t (\mathcal L_t))$ that $\Phi_t(\mathcal L_t)$ is
in strong generic position for each $t$ in a neighbourhood $W$ of $o$.
Hence the topological type of $\pi(\Phi_t (\mathcal L_t))$ equals the topological type of $\pi(\Phi_o (\mathcal L_o))$ for $t\in W$. 
 Moreover,  the multiplicity of the reduced curve $\pi(\Phi_t(\mathcal L_t))$  is constant near $o$. Therefore $\pi_T (\Phi(\mathcal L))_o$ is reduced.
$\square$  
\end{proof}

\begin{e-definition}\label{microlocal}\em 
Let $Y$ be a plane curve with an irreducible tangent cone. 
Let $L$ be the conormal of $Y$.
We say that an [equisingular] microlocal deformation $\mathcal{Y}$ of $Y$ over $T$ 
is an  [\em equisingular\/\em ] \em microlocal versal deformation of  $Y$ \em if 
{\em $($i$)$}  
for each [equisingular] microlocal deformation $Y\hookrightarrow \mathcal Z \to (R,o)$, there is $\psi:R\rightarrow T$ such that 
 $\psi^{\ast}\mathcal{Y}$ is microlocally equivalent to $\mathcal Z$; {\em $($iii$)$} 
the derivative of $\psi$ at  $o$ only depends on $\mathcal Z$.
\end{e-definition}

\begin{e-definition}\em 
An [equisingular] deformation $\mathcal L$ of a Legendrian curve $L$  over $T$ 
is an   [\em equisingular\em\/] \em versal deformation of  $L$ \em if {\em $($i$)$} for each [equisingular] deformation 
$\mathcal K$ of $L$ over $(R,o)$ there is $\psi : R\rightarrow T$ such that 
 $\psi^{\ast}\mathcal{L}$ is isomorphic  to $\mathcal K$; {\em $($ii$)$} 
The derivative of $\psi$ at  $o$ only depends on $\mathcal K$.
\end{e-definition}

\noindent
Let $k,n$ be integers such that $n> k>1$ and $(k,n)=1$. Set $Y=\{y^k-x^n=0\}$.
Let $c$ be the conductor of the semi-group of  $Y$. Let $\mathcal L$ be an equisingular deformation of the conormal $L$ of $Y$.

\begin{remark}\label{lemaone}\em
An equisingular deformation $\mathcal Y$ of $Y$ admits a parametrization of the type $x=s^k$, 
$y=s^n+$ $\sum_{r\ge n+1}a_rs^r$, where $a_r\in\mathcal O_{T,o}$ and $a_r(o)=0$.
Hence the conormal $\mathcal K$ of $\mathcal Y$ admits a parametrization 
  $\sigma : (\mathbb C,0)\times (T,o) \to \mathcal K$ of type (\ref{expfamilia2}).
For each $t$, $s\mapsto \sigma(s,t)$ defines a parametrization of $\mathcal K_t$, 
hence $\mathcal K_t$ is the conormal of $\mathcal Y_t$. Therefore  $\mathcal Y$ is a microlocal deformation of $Y$. 
The parametrization $\sigma$ defines a valuation $w$ of the ring $\mathbb C[[x,y,p,t]]$, 
the one that associates to $f\in\mathbb C[[x,y,p,t]]$ the multiplicity of the zero of $f\circ\sigma$ has an element of $\mathbb C[[t]][t^{-1}][[s]]$.
\end{remark}


\noindent
\begin{lemma}\label{cleaning}
Let $i,j,l$ be non negative integers such that $ki+nj+(n-k)l>kn$. There are non negative integers $a,b$ such that $w(x^iy^jp^l)=w(x^ay^b)$, 
$x^iy^jp^l\equiv (n/k)^l x^ay^b$ mod $(t)+I_{\mathcal L}$ and $w(x^iy^jp^l-(n/k)^l x^ay^b)>w(x^ay^b)$.
Moreover, $p^k\equiv (n/k)^k x^{n-k}$ mod $(t)+I_{\mathcal L}$.
\end{lemma}

\begin{proof}
By remark \ref{lemaone}, $xp\equiv (n/k)y$, $x^{i+1} y^j p^{l+1}\equiv (n/k) x^{i} y^{j+1} p^{l}$ and $p^k\equiv (n/k)^k x^{n-k}$ mod $(t)+I_{\mathcal L}$.  
Hence we can assume that $l<k$.
If $i=0$, then $j+l>k$ and  $y^j p^l\equiv (n/k)^l x^{n-l} y^{j+l-k}$
mod $(t)+I_{\mathcal L}$.  
\end{proof}

\noindent
Given  an analytic set $Z$ of $(T,o)$, let $I_{Z}$[$\widehat{I}_{Z}$] 
be the ideal of elements of $\mathcal O_{T,o}$[$\widehat{\mathcal O}_{T,o}$] that vanish on $Z$.

\begin{lemma}\label{inteirosformais}
Given $u\in \mathbb C\{ x,y,p,t\}$ such that $w(u)\geq c$, there is $v\in  \mathbb C[[ x,y,p,t]] $ 
such that $u-v\in \widehat{I}_{\mathcal{L}}$.
\end{lemma}

\begin{proof}
There are non negative integers $a,b$ and $\xi\in\mathbb C\{t\}$ 
such that $w(u)=w(x^a y^b)$ and $w(u-\xi x^a y^b)> w(u)$. 
We iterate the procedure, producing in this way $v\in\mathbb C[[x,y,t]]$ such that $w(u-v)=\infty$.$\square$
\end{proof}

\begin{lemma}\label{contactbuild} \em (cf. \cite{ARNE}) \em
Given non negative integers $a,b$ and 
$\lambda\in\mathbb C$, there are a contact transformation $\varphi$ of type $(\ref{grupojota})$ and $\varepsilon\in \mathbb C\{x,y,p\}$ 
such that  $\alpha=\lambda y^{a}p^{b}$,
$\beta=\lambda by^{a} p^{b+1}/(b+1)+\varepsilon$
and $w(\varepsilon)\geq w(y^{2a-1}p^{2b+2})$. 
\end{lemma}



\begin{theorem}\label{versaldef} 
The function
$
F(x,y,t)=y^k-x^n+\sum_{(i,j)\in B}t_{i,j}~ x^iy^j
$
defines an equisingular versal deformation $\mathcal F$ of $\{y^k-x^n=0\}$ over $\mathbb C^{ B}$, 
where $B=\{ (i,j) ~: ~ki+nj > kn,~i\le n-2, ~j\le k-2 \}$.
\end{theorem} 
\begin{proof}Its a corollary of theorem II.1.16 of \cite{GP5}.\end{proof}$\square$

\begin{theorem}\label{microykxn}
If $n>2k$, the function 
$
G(x,y,t)=y^k-x^n+{\sum_{(i,j)\in C}}t_{i,j}~ x^iy^j
$
defines an equisingular microlocal versal deformation $\mathcal G$ of $\{y^k-x^n=0\}$ over $\mathbb C^{ C}$,
where $C=\{ (i,j)\in B ~: ~i+j\leq n-2 \}$.
\end{theorem}


\begin{proof}
We can assume that $R=(\mathbb C^a,0)$ and $\mathcal O_{\mathbb C^{a},0}=\mathbb C\{r\}$. 
Let $H\in\mathbb C\{x,y,r\}$ be a generator of the ideal $I_{\mathcal Z}$. 
By theorem \ref{versaldef} we can assume that there are
 $\xi_{i,j}\in\mathbb C\{r\}$, 
$(i,j)\in B$,  such that $\xi_{i,j}\in (r)$ and
$
H(x,y,r)=y^k-x^n+\sum_{(i,j)\in B}\xi_{i,j}x^iy^j
$.
Moreover,
the functions $\xi_{i,j}$ are uniquely determined mod $(r^2)$. 
Let $\mathcal K$ be the conormal of $\mathcal Z$.
Assume that there are $\psi:\mathbb C^a\rightarrow\mathbb C^{ C}$, $\psi=(\psi^{i,j})$, and $\Phi\in \mathcal J'_{\mathbb C^a}$ such that
$\Phi^{-1} (\mathcal K)$ equals the conormal of $\psi^*\mathcal G$. 
There are $\alpha,\beta,\gamma\in \mathbb C\{x,y,p,r\}$
such that $\Phi$ equals (\ref{grupojota}).
Hence there is $\varepsilon\in\mathbb (r)$  such that
$
H(x+\alpha,y+\beta,r)=(1+\varepsilon(x,y,r))G(x,y,\psi(r))\ \ \textrm{mod}\ I_{\Phi^{-1}(\mathcal K)}.
$
Notice that 
$$
H(x+\alpha,y+\beta,r)\equiv H(x,y,r)+ky^{k-1}\beta-nx^{n-1}\alpha \ \ \textrm{mod}\  (r)^2.
$$
By lemma \ref{cleaning}, $x^{n-1}\alpha$, $y^{k-1}\beta$ are congruent modulo  $(r)^2+\widehat{I}_{\Phi^{-1}(\mathcal K)}$
with elements of  $(x,y)^{n-1}\mathbb C[[x,y,r]]$.
Hence $\psi^{i,j}\equiv\xi_{i,j}$ mod $(r)^2$ for each $(i,j)\in C$.\hphantom{x}

\noindent 
By theorem \ref{versaldef}, 
it is enough to show that there is 
$\psi:\mathbb C^{ B}\rightarrow\mathbb C^{ C}$ such that $\psi^*\mathcal G$ is microlocally equivalent to $\mathcal F$.
Set $N=\# (B\setminus C)$. Let us order the pairs $(i,j)\in B\setminus C$ by the value of $ki+nj$. 
Let $B_{l}$ be the set $B$ minus the set of the $l$ smaller ordered pairs of $B\setminus C$. 
Set 
$H_{l}(x,y,t_{B_l})= y^k-x^n+\sum_{(i,j)\in B_{l}}t_{i,j}x^i y^j$. 
Let $\mathcal L_l$ be the conormal of $\mathcal H_l=\{H_l=0\}$.
Let $0\leq l\leq N-1$. Its enough to show that there is $\psi_l:\mathbb C^{ B_l}\to \mathbb C^{ B_{l+1}}$ such that
$\psi_l^*\mathcal H_{l+1}$ is microlocally equivalent to $\mathcal H_l$.
Let $(a,b)$ be the smallest element of $B_l\setminus C$. 
Let $v=w_l(x^a y^b)$. 
Let $\Phi_l$ be the contact transformation of the type (\ref{grupojota})
associated by theorem \ref{ABC} to $\alpha=\lambda y^{a+b-(n-1)}p^{n-1-a}$ and $\beta_0=0$.
Let $w_l$ be the valuation associated to $\mathcal \Phi_l^{-1}(\mathcal L_l)$.
By lemma \ref{contactbuild}, $\beta=\lambda (n-1-a)y^{a+b-(n-1)}p^{n-a}/(n-a)+\mu$, 
with $w_l(\mu)>w_l(y^{a+b-(n-1)}p^{n-a})$. 
There is $\delta\in \mathbb C\{x,y,p,t_{B_l}\}$ such that $w_l(\delta)> v$ and 
\begin{displaymath}
H_l(x+\alpha,y+\beta,t_{B_l})\equiv 
H_{l+1}(x,y,t_{B_{l+1}})+t_{a,b}x^{a} y^{b}+ky^{k-1}\beta-nx^{n-1}\alpha+\delta\ \ \textrm{mod} \ I_{\Phi_l^{-1}(\mathcal L_l)}.
\end{displaymath}
Set $\lambda=t_{a,b}(n-a)(n/k)^{n-1-a}/n$. 
By lemma \ref{inteirosformais}, there is $\widehat{\delta}\in \mathbb C[[x,y,t_{B_l}]]$ such that $w_l(\widehat{\delta})>v$ and  
$H_{l}(x+\alpha,y+\beta,t_{B_l})\equiv 
H_{l+1}(x,y,t_{B_{l+1}})+\widehat{\delta}$ 
mod $\widehat{I}_{\Phi_l^{-1}(\mathcal L_l)}$. 
Since $\pi_{B_{l}}(\Phi_l^{-1}(\mathcal L_l))$ is a convergent hypersurface, 
we can assume that $\widehat{\delta}\in \mathbb C\{x,y,t_{B_l}\}$. 
By the proof of theorem \ref{versaldef}, there are $\alpha,\beta,\varepsilon\in\mathbb C\{x,y\}[[t_{B_l}]]$ and 
$\psi^{i,j}\in \mathbb C[[t_{B_l}]]$, $(i,j)\in B$, such that $\alpha,\beta,\varepsilon\in (t_{B_l})$ and
\begin{displaymath}
(H_{l+1}+\widehat{\delta})(x+\alpha,y+\beta,t_{B_l})=(1+\varepsilon )F(x,y,\psi_{l}(t_{B_l})),
\end{displaymath}
where $\psi_l=({\psi^{i,j}})$.
Since $y^k-x^n$ is quasi-homogeneous, we can assume that $\varepsilon$ vanishes.
Set $\alpha=\sum_q\alpha_q$, $\beta=\sum_q\beta_q$,  ${\psi^{i,j}}=\sum_q\psi^{i,j}_q$,
where $\alpha_q,\beta_q,\psi^{i,j}_q$ are homogeneous functions of degree $q$ on the variables $t_{i,j}$, $(i,j)\in B_l$. 
For each $q$ there is an homogeneous part $g_q$ of degree $q$  of an element of the ideal $(\partial_x H_{l+1},\partial_y H_{l+1})$
such that $\alpha_q nx^{n-1}   +\beta_q   ky^{k-1} = g_q$.
Moreover, each $g_q$ depends on the choices of $\alpha_p,\beta_p$, $p<q$.
Notice that $\alpha_0,\beta_0,g_0=0$, $g_1$ is the linear part of $\widehat{\delta}$, 
$\psi^{i,j}_1=t_{i,j}$ if $(i,j)\in B_{l+1}$ and $\psi^{i,j}_1=0$ otherwise.
We choose  $\alpha_1,\beta_1$ such that $\beta_1ky^{k-1}$ is the rest of the division of $g_1$ 
by  $nx^{n-1}$. 
Since $w_l(g_1)>v$, then $w_l(\alpha_1x^{n-1}),w_l(\beta_1y^{k-1}),w_l(g_2)>v$ and  $\psi^{i,j}_2=0$ if $ki+nj\le v$.
Moreover, we can iterate the procedure. We show in this way that $F(x,y,\psi_{l}(t_{B_l}))=H_{l+1}(x,y,\psi_{l}(t_{B_l}))$.
Since our choices of $\alpha_q,\beta_q$ are the ones of lemma 1 of \cite{KaS}  and of theorem II.D.2 of \cite{GR}, the functions
$\alpha,\beta, \psi^{i,j}$, $(i,j)\in B$,  converge.
$\square$
\end{proof}

\begin{corollary}
Each irreducible Legendrian curve $L$ contained in a smooth surface admits an equisingular versal deformation.
This deformation is trivial if and only if $L$ is isomorphic to the conormal of a curve defined by one of the functions
$y^2-x^{2n+1}$, $n\ge 1$, $y^3-x^7$, $y^3-x^8$.
\end{corollary}
\begin{proof}
By theorem 8.3 of \cite{SKKO}, there are integers $k,n$ such that $n>2k>1$, $(k,n)=1$ 
and the Legendrian curve is isomorphic to the conormal of $\{y^k=x^n\}$. The corollary  follows from results \ref{lemaone}, 
 \ref{lematwo} and \ref{microykxn}.
$\square$
\end{proof}

\noindent
By the canonical properties of versal deformations (cf. theorem II.1.15 of \cite{GP5}), 
the previous result still holds if $L$ is the conormal of a semi-quasi-homogeneous branch.

\begin{example}Set
$f_0(x,y)=y^4-x^{11}$, $f_1(x,y)=y^4-x^{11}+x^6y^2$, $f_2(x,y)=y^4-x^{11}+x^7y^2$. 
Let $L_i$ be the conormal of $\{f_i=0\}$, $0\le i\le 2$.
A Legendrian curve equisingular to $L_0$ is isomorphic to one and only one of the curves $L_0,L_1,L_2$.
\em
The equisingular microlocal versal deformation of $y^4-x^{11}$ equals
$G(x,y,t_2,t_6)=y^4-x^{11}+t_2x^6y^2+t_6x^7y^2$.  
Moreover,
\begin{equation}\label{ddeeffppp}
G(\lambda^4 x,\lambda^{11} y,t_2,t_6)=\lambda^{44}G(x,y,\lambda^2t_2,\lambda^6t_6), \qquad \textrm{for each $\lambda\in \mathbb C^*$.}
\end{equation}
If $t_2\not=0$, it follows from (\ref{ddeeffppp}) that we can assume $t_2=1$. 
Moreover, there are $\varepsilon,\delta\in (x,y)$ such that 
$G(x-2t_6 x^2, y- 11 t_6xy/2,1,t_6)= (1-22t_6 x+\varepsilon )(y^4-x^{11}+ x^6y^2+\delta)$ and $w(\delta)\geq 52 $.
By theorem \ref{microykxn} we can assume that $\delta =0$.
If $t_2=0$ it follows from (\ref{ddeeffppp}) that we can assume $t_6=0$ or $t_6=1$. 

\noindent
Since $L_0$ is contained in the surface $11y-4xp=0$ 
and $L_1$, $L_2$ are not contained in a smooth surface, 
$L_0$ cannot be isomorphic to $L_1$ or $L_2$. Let $\varphi\in\mathcal J$.
There is $\varepsilon\in\mathbb C\{x,y\}$ such that 
$f_1(x+\alpha,y+\beta)\equiv f_1(x,y)+\varepsilon(x,y)$ mod $I_{\varphi(L_1)}$.
If $f_1+\varepsilon=f_2$, $w(\varepsilon )\ge 46$. Hence $w(\varepsilon )\ge 47$. Therefore $f_1+\varepsilon\not=f_2$.
By theorem \ref{inv}, $\varphi (L_1)\not=L_2$ for each $\varphi\in\mathcal G$. 
\end{example}

\bigskip\noindent
Orlando Neto and Jo\~ao Cabral, 

\noindent
Microlocal versal deformations of the plane curves $y^k = x^n$, 

\noindent
C. R. Acad. Sci. Paris, Ser. I 347 (2009), pp 1409 $-$ 1414.

\bigskip\noindent
http://ptmat.fc.ul.pt/$\sim$orlando/publications.html
\end{document}